# Estimation Technique for a Contact Point Between two Materials in a Stationary Heat Transfer Problem

Guillermo Federico Umbricht[1,2*], Diana Rubio[1], Domingo Alberto Tarzia[3,4]

[1] Centro de Matemática Aplicada, Escuela de Ciencia y Tecnología, Universidad Nacional de San Martín. 25 de Mayo y Francia, San Martín, Buenos Aires, B1650, Argentina
[2] Instituto de Ciencias e Instituto del Desarrollo Humano, Universidad Nacional de Gral. Sarmiento. Juan María Gutiérrez 1150, Los Polvorines, Buenos Aires, B1613, Argentina
[3] Consejo Nacional de Investigaciones Científicas y Técnicas (CONICET), Godoy Cruz 2290, Ciudad Autónoma de Buenos Aires,, C1425FQB, Argentina
[4] Departamento de Matemática, Facultad de Ciencias Empresariales, Universidad Austral, Paraguay 1950, Rosario, Santa Fe, S2000FZF, Argentina

*Keywords:*
Elasticity analysis, heat transfer, interface problem, mathematical modeling, numerical simulation, parameter estimation.

**ABSTRACT**

An inverse problem for a stationary heat transfer process is studied for a totally isolated bar on its lateral surface, of negligible diameter, made up of two consecutive sections of different, isotropic and homogeneous materials. At the left boundary, a Dirichlet type condition is imposed that represents a constant temperature source while a Robin type condition that models the heat dissipation by convection is considered at the right one. Many articles in the literature focus on thermal and stress analysis at the interface but no one is dedicated to the estimation of the contact point location between the two materials. In this work, it is assumed that the interface position is unknown. A technique to determine it from a unique noisy flow measurement at the right boundary is introduced. Necessary and sufficient conditions are derived in order to obtain the estimation of the interface point from a heat flux measured at the right boundary. Numerical solutions are obtained together with an expression for the estimation error. Moreover, an elasticity analysis is included to study the influence of data errors. The results show that our approach is useful for determining the location of the materials interface.

## 1. INTRODUCTION

Heat transfer problems in multilayer or solid-solid interface materials have been studied in recent years due to the multiple and different applications that have been found in science and engineering [1].

These problems have direct applications in different industries, among the most important, the metallurgical [2], the technological and electronic [3], the automotive [4], aerospace [5] and aviation [6].

The advancement of technology requires materials with particular thermal, electrical, magnetic, acoustic and optical properties due to which the interface properties of different combinations of materials have been studied. Among the most studied are the following pairs Al-Cu [7], Cu-Cu [7, 8], Cu-Al [9], Pb-Al [2]. Other materials are considered ( [6], [8], [10- 15]). These works focus on tension [7], the adhesion [16], thermal resistance [8, 17], the corrosion [6], electrical conductivity [2] and thermal conductivity [12]. Many properties have been studied at the interface of different materials but there are no techniques in the literature that allows to determine the contact point of perfectly joined materials.

In this work, we consider the stationary problem of heat transfer of a bar of negligible diameter and known length, totally isolated on its lateral surface, composed of two different, isotropic and homogeneous materials. It is assumed that the temperature on the left boundary of the bar is controlled by a thermal source that is maintained at a constant temperature, the thermal resistance at the interface is neglected and the right edge is left free, giving rise to the phenomenon of convection.

From the analytical solution of the forward problem, the interface position is estimated by means of a unique heat flux measurement at the right edge of the bar. Necessary and sufficient conditions are derived in order to obtain the estimation of the interface point.

In order to study the local influence of the flow measurement in the estimation, an elasticity analysis is performed.

Numerical examples considering different situations are included and commented to illustrate the results for the estimation technique proposed here.

## 2. MATHEMATICAL MODEL

In this section, a mathematical model for the interface problem is stated and an analytical expression for the solution is found. Furthermore, the analytical solution is shown to be consistent with the one for a homogenous bar, and temperature profiles are shown for particular cases.

## 2.1 Model

The steady-state heat transfer problem [18, 19] for a bar of known length, negligible diameter, fully insulated on its lateral surface, composed of two consecutive, different, isotropic and homogeneous materials, can be modeled by

$$u'' = 0, \quad 0 < x < l, \tag{1}$$

$$u'' = 0, \quad l < x < L, \tag{2}$$

where $u$ (°C) represents the temperature of the bar, $L$ (m) the length of the bar and $l$ (m) the location of the contact point. It is assumed that the section of the bar made with a material $A$ has length $l$ and the section of the bar made with the material $B$ has length $L-l$.

Thermal resistance is assumed to be neglected; hence, conditions of equal temperature and flow are imposed at the interface [19], that is

$$u(l^+) = u(l^-), \tag{3}$$

$$\kappa_B u'(l^+) = \kappa_A u'(l^-), \tag{4}$$

where $\kappa_A, \kappa_B$ (Wm$^{-1}$°C$^{-1}$) represent the thermal conductivities of the materials $A$ and $B$, respectively.

Boundary conditions are given by a source with constant temperature on the left and convective condition on the right [19]

$$u(x) = F, \quad x = 0, \tag{5}$$

$$\kappa_B u'(x) = -h(u(x) - T_a), \quad x = L, \tag{6}$$

where $F$ (°C) represents the source at constant temperature, $T_a$ (°C) the room temperature and $h$ (Wm$^{-2}$°C$^{-1}$) the convection heat transfer coefficient. In this work it is assumed that $F > T_a$.

Fig. 1 outlines the problem of interest.

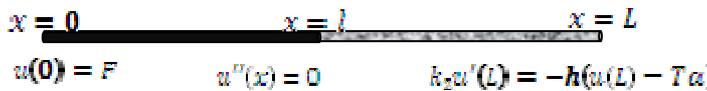

**Figure 1.** Scheme for the mathematical model

## 2.2 Analytical solution of the Forward Problem

The problem to be solved is described by equations (1)-(6).

From expressions (1)-(2) it follows that

$$u(x) = a + bx, \quad 0 < x < l, \tag{7}$$

$$u(x) = c + dx, \quad l < x < L, \tag{8}$$

where $a, b, c$ and $d$ are constants to be determined from the conditions given in (3)-(6)-as follows

$$a + bl = c + dl, \tag{9}$$

$$\kappa_B d = \kappa_A b, \tag{10}$$

$$a = F, \tag{11}$$

$$\kappa_B d = -h(c + dL - T_a). \tag{12}$$

Equations (9)-(12) leads to

$$a = F, \tag{13}$$

$$b = \frac{\kappa_B h(T_a - F)}{\varsigma}, \tag{14}$$

$$c = F + \frac{lh(\kappa_B - \kappa_A)(T_a - F)}{\varsigma}, \tag{15}$$

$$d = \frac{\kappa_A h(T_a - F)}{\varsigma}, \tag{16}$$

where

$$\varsigma = \kappa_A \kappa_B + \kappa_A hL + (\kappa_B - \kappa_A)hl. \tag{17}$$

Replacing expressions (13)-(17) in (7)-(8) the analytical solution of the problem of interest is obtained, which is given by

$$u(x) = F + \frac{\kappa_B h(T_a - F)}{\varsigma} x, \quad 0 \leq x \leq l, \tag{18}$$

$$u(x) = F + \frac{h(T_a - F)(l(\kappa_B - \kappa_A) + \kappa_A x)}{\varsigma}, \quad l < x \leq L.$$

## 2.3 Model consistency

In the particular case where the bar is made of only one material ($\kappa_A = \kappa_B = \kappa$) the solutions given in equations (18) are reduced to

$$u(x) = F + \frac{\kappa h(T_a - F)}{\varsigma} x, \quad 0 \leq x \leq l, \tag{19}$$

$$u(x) = F + \frac{h(T_a - F)\kappa}{\varsigma} x, \quad l < x \leq L, \tag{20}$$

or, equivalently,

$$u(x) = F + \frac{\kappa h(T_a - F)}{\varsigma} x, \quad 0 \leq x \leq L. \tag{21}$$

The solution given in (21) coincides with the solution to the stationary heat transfer problem with the same boundary conditions imposed for a homogeneous bar [19].

## 2.4 Examples

In order to illustrate the temperature profiles for the problem described by (1)-(6), few examples are considered, where the materials involved and the contact points are different for each case. In all of them, it is assumed a bar of $L=1$ m with $F=100\,°C$ $T_a = 25\,°C$ and $h=10$ Wm$^{-2}$°C$^{-1}$. The average thermal conductivity values used for the examples were obtained from [4] and they are included in Tab. 1.

**Table 1.** Thermal conductivity of different materials

| Material | symbol | $\kappa$ (Wm$^{-1}$°C$^{-1}$) |
|---|---|---|
| Aluminium | Al | 204 |
| Cupper | Cu | 386 |
| Iron | Fe | 73 |
| Silver | Ag | 419 |
| Lead | Pb | 35 |
| Magnesium | Mg | 156 |

In Figure 2, four plots are shown corresponding to Al-Cu ($l = 0.3m$ and $l = 0.5m$) and Cu-Al ($l = 0.3m$ and $l = 0.5m$) for a bar of length $L = 1m$. When the interface is in the middle of the bar ($l = 0.5m$) the relative location of the two materials (Al-Cu or Cu-Al) makes no difference on the temperature value reached on the right boundary, which, in this case, achieves 97.25°C. This observation is valid for any pair of materials whenever $l = L/2$.

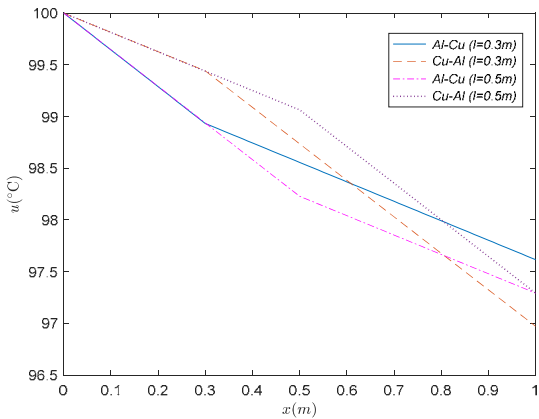

**Figure 2.** Temperature for different materials

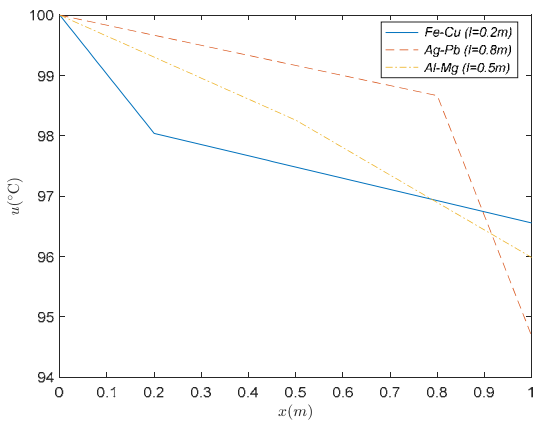

**Figure 3.** Temperature for different materials

In Fig. 3 you can see the differences in the temperature profiles in three particular cases: $\kappa_A < \kappa_B$ (Fe-Cu), $\kappa_A > \kappa_B$ (Ag-Pb) and similar values for $\kappa_A$ y $\kappa_B$ (Al-Mg). This results are consistent with the temperature gradients given in equations (14) and (16).

## 3. ESTIMATION OF THE POINT OF CONTACT

In this section, the main objective of this work is developed. Under the conditions given by equations (1)-(6) the contact point position ($x = l$) is estimated using a noisy flow data taken at $x = L$. In addition, the necessary and sufficient conditions for the estimation are provided and a bound is given for the estimation error.

### 3.1 Estimation of the interface location

By definition, the thermal flux $q$ at $x = L$ is given by the following expression [19]

$$q = -\kappa_B u'(x), \quad x = L. \tag{22}$$

From equation (18) along with (22) the following expression of $q$ is obtained for the particular problem addressed in this work, by the following expression

$$q = \frac{\kappa_B \kappa_A h (F - T_a)}{\kappa_A \kappa_B + \kappa_A h L + (\kappa_B - \kappa_A) h l} \tag{23}$$

then, the parameter $l$ is obtained from equation (23) as

$$l = \frac{\kappa_A \kappa_B}{\kappa_B - \kappa_A} \left( \frac{F - T_a}{q} - \frac{1}{h} - \frac{L}{\kappa_B} \right). \tag{24}$$

Therefore, the estimate $\hat{l}$ of the contact point as a function of the flow measurement $\hat{q}$ at $x=L$ is expressed as

$$\hat{l} = \frac{\kappa_A \kappa_B}{\kappa_B - \kappa_A} \left( \frac{F - T_a}{\hat{q}} - \frac{1}{h} - \frac{L}{\kappa_B} \right), \tag{25}$$

where it is assumed that

$$|q - \hat{q}| \leq \varepsilon, \tag{26}$$

where $\varepsilon$ represents the noise level in the data, which, in practice, it is easily determined from the error in the measurement instruments used.

Note that the estimation of $l$ only depends on the flow measurement and the parameters of the problem.

### 3.2 Necessary and sufficient conditions

There are necessary and sufficient conditions for the estimation of the interface point. Regardless of the noise in the measurement of $q$, the value of $\hat{l}$ must satisfy

$$0 < \hat{l} < L, \tag{27}$$

or, equivalently, from equation (25)

$$0 < \frac{\kappa_A \kappa_B}{\kappa_B - \kappa_A}\left(\frac{F-T_a}{\hat{q}} - \frac{1}{h} - \frac{L}{\kappa_B}\right) < L. \quad (28)$$

This expression imposes conditions on the values of $\hat{q}$. To find these conditions, two cases are studied.

Case 1 ($\kappa_A < \kappa_B$): Since $0 < \kappa_B - \kappa_A$, from (28) it results

$$\frac{\kappa_A h(F-T_a)}{Lh + \kappa_A} < \hat{q} < \frac{\kappa_B h(F-T_a)}{Lh + \kappa_B}, \quad (29)$$

Case 2 ($\kappa_A < \kappa_B$): Now, $0 > \kappa_B - \kappa_A$, and from (28) it results

$$\frac{\kappa_B h(F-T_a)}{Lh + \kappa_B} < \hat{q} < \frac{\kappa_A h(F-T_a)}{Lh + \kappa_A}. \quad (30)$$

From equations (29) - (30) arise the following necessary condition for the estimation of the contact point given by

$$q_m = \frac{\kappa_m h(F-T_a)}{Lh + \kappa_m} < \hat{q} < \frac{\kappa_M h(F-T_a)}{Lh + \kappa_M} = q_M, \quad (31)$$

where

$$\kappa_m = \min\{\kappa_A, \kappa_B\}, \kappa_M = \max\{\kappa_A, \kappa_B\}. \quad (32)$$

**Remark:** Inequalities (31) give us a necessary and sufficient conditions for the estimation of the interface point $l$ ($0 > l > L$).

### 3.3 Error estimate

Assuming that the measured flow $\hat{q}$ satisfies the conditions given by the expressions (26) and (31)-(32), an analytical bound for the error in the estimation of the contact point follows for $K = K(\hat{q})$ such that

$$\left|l - \hat{l}\right| \le K. \quad (33)$$

By replacing the expressions (24)-(25) in (33), the left side hand becomes

$$\left|l - \hat{l}\right| = \frac{\kappa_A \kappa_B}{|\kappa_B - \kappa_A|} \frac{F-T_a}{q\hat{q}} |q - \hat{q}|. \quad (34)$$

From here and (26), it follows that a bound can be chosen according to

$$K = \frac{\kappa_A \kappa_B}{|\kappa_B - \kappa_A|} \frac{F - Ta}{q\hat{q}} \varepsilon. \quad (35)$$

Note that the bound $K$ could becomes very large when the materials have similar thermal conductivities, i.e. $\kappa_A - \kappa_B \approx 0$.

## 4. LOCAL DEPENDENCE OF THE PARAMETER WITH RESPECT TO THE DATA

Equation (25) indicates that the estimated value for $l$ depends on the parameters of the problem and on the measured flow $\hat{q}$. There are some tools that helps to study the influence of data $\hat{q}$ on the estimated parameter $\hat{l}$. Some of the most useful ones are the sensitivity [20] and the elasticity analysis [21]. In this work, the latter one is applied.

### 4.1 Elasticity

This technique is widely used in economics. It provides the percentage error in the estimated parameter for 1 % error in a measurement value. It is defined by

$$E(q) = \frac{q}{l}\frac{\partial l}{\partial q}, \quad (36)$$

hence, the expression (24) yields

$$E(q) = \frac{(F-T_a)h\kappa_B}{q(\kappa_B + Lh) - h\kappa_B(F-T_a)}. \quad (37)$$

### 4.2 Elasticity function analysis

The elasticity function given by expression (37) has particularities that deserve to be highlighted.

#### 4.2.1 Vertical Asymptote
The vertical asymptote for the function (37) is given by the expression

$$q = \frac{h\kappa_B(F-T_a)}{\kappa_B + Lh}. \quad (38)$$

Note that this can be interpreted in two different cases. From equations (29)-(32), if $\kappa_A < \kappa_B$, then the elasticity has a vertical asymptote at $q = q_M$. On the other hand, if $\kappa_A > \kappa_B$ then the vertical asymptote is at $q = q_m$.

#### 4.2.2 Sign
Under the conditions studied here, the numerator of the elasticity function is always strictly positive, then the sign of the elasticity function depends on the denominator.

From expressions (37), taking into account the expressions given by (29) and (30), it follows that

$$E(q) < 0 \Leftrightarrow q < \frac{h\kappa_B(F-T_a)}{\kappa_B + Lh} \Leftrightarrow \kappa_A < \kappa_B, \quad (39)$$

$$E(q) > 0 \Leftrightarrow q > \frac{h\kappa_B(F-T_a)}{\kappa_B + Lh} \Leftrightarrow \kappa_A > \kappa_B. \quad (40)$$

In other word, $\kappa_A < \kappa_B$ yield to a negative elasticity function in the interval $[q_m, q_M]$, and conversely. Otherwise, if $\kappa_A > \kappa_B$ the function turns out to be positive in the same interval, and reciprocally.

### 4.2.3 Monotony

Another important observation is that the elasticity function turns out to be decreasing. This fact can be easily seen by differentiating the expression (37) given by

$$\frac{\partial E(q)}{\partial q} = -\frac{(F-T_a)h\kappa_B(\kappa_B+Lh)}{\left[q(\kappa_B+Lh)-h\kappa_B(F-T_a)\right]^2} < 0 \quad (41)$$

Note that (41) implies that $E(q)$ is a strictly decreasing function.

## 5. NUMERICAL EXAMPLES

Numerical examples corresponding to the three cases: $\kappa_A < \kappa_B$, $\kappa_A > \kappa_B$ and similar values for $\kappa_A$ and $\kappa_B$, are included.

For all these numerical examples, the following values are imposed, $L=$ 10 m, $F=$100 °C, $T_a = 25$ °C, $h=$10 Wm$^{-2}$°C$^{-1}$. The average thermal conductivities were obtained from [4] and are included in the Tab. 1.

For this work, the value of $q$ is calculated from the forward problem (1) - (6) together with equation (22). Noise is added to numerically simulate the experimental measurement $\hat{q}$. Then, from equation (25), the estimate value for the contact point $\hat{l}$ is found.

### 5.1 Example 1

In this example the case $\kappa_A < \kappa_B$ is considered. In particular, a Fe-Cu bar is assumed to estimate the interface position $l$ ($\kappa_A = 73$Wm$^{-1}$ $^0$C$^{-1}$ < $\kappa_B = 386$Wm$^{-1}$ $^0$C$^{-1}$).

For $q = 440.299$ W.m$^{-2}$, the necessary and sufficient conditions (31) are verified since $q_m = 316.474$ W.m$^{-2}$ and $q_M = 595.679$ W. m$^{-2}$ and from equation (24) it follows that $l = 4m$.

Since usually data are noisy due to error measurements, we calculate the errors in the estimation of the location $l$ for several heat flux values that satisfy the necessary and sufficient conditions (31), that is, for $q_m < q < q_M$.

The Tab. 2 contains noisy value for flow measurements in a range around the actual value $q$, along the level noise $\varepsilon$, the corresponding estimates $\hat{l}$ for $l$ and the estimated error bound $K$ defined in (35).

**Table 2.** Estimate of $l$ for Example 1.

| $\hat{q}$ (W. m$^{-2}$) | $\hat{l}$ (m) | $\varepsilon$ (W. m$^{-2}$) | $K$(m) |
|---|---|---|---|
| 436 | 4.151 | 4.299 | 0.151 |
| 437 | 4.115 | 3.299 | 0.115 |
| 438 | 4.080 | 2.299 | 0.080 |
| 439 | 4.045 | 1.299 | 0.045 |
| 440.299 | 4.000 | 0.000 | 0.000 |
| 441 | 3.975 | 0.701 | 0.024 |
| 442 | 3.941 | 1.701 | 0.059 |
| 443 | 3.906 | 2.701 | 0.093 |
| 444 | 3.872 | 3.701 | 0.127 |
| 445 | 3.838 | 4.701 | 0.162 |

In Fig. 4 the elasticity is plotted for Example 1. Due to the fact $\kappa_A < \kappa_B$, the elasticity is a negative function as shown in (39), strictly decreasing (see (41)) with a vertical asymptote at $q = q_M$ (see (38)). This function indicates that a measurement error of 1 % in the flow value $\hat{q}$ translates into an error of around 4 % in the estimation value $\hat{l}$ of the contact point.

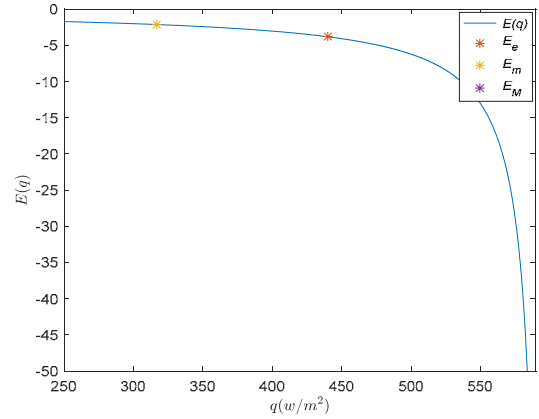

**Figure 4.** Elasticity for Example 1

### 5.2 Example 2

In this example we consider case $\kappa_A > \kappa_B$. In particular, a Ag-Pb bar is assumed to estimate the interface length $l$ ($\kappa_A = 419$Wm$^{-1}$ $^0$C$^{-1}$ < $\kappa_B = 35$Wm$^{-1}$ $^0$C$^{-1}$).

For $q = 266.927$ W.m$^{-2}$, the necessary and sufficient condition (31) are verified since that $q_m = $ 194.444 W.m$^{-2}$ and $q_M = 605.491$ W.m$^{-2}$ and from equation (24) it follows that $l = 4m$.

The Tab. 3 contains the noisy value for flow measurements in a range around the actual value $q$, along the level noise $\varepsilon$ with the corresponding estimates $\hat{l}$ for $l$ and the estimated error bound $K$ defined in (35).

In Fig. 5 the elasticity is plotted for example 2. Due to the fact $\kappa_A > \kappa_B$, the elasticity is a positive function (see (40)), strictly decreasing (see (41)) having a vertical asymptote at $q = q_m$ (see (38)).

**Table 3.** Estimate of $l$ for Example 2

| $\hat{q}$ (W. m$^{-2}$) | $\hat{l}$ (m) | $\varepsilon$ (W. m$^{-2}$) | $K$(m) |
|---|---|---|---|
| 263 | 3.839 | 3.927 | 0.160 |
| 264 | 3.881 | 2.927 | 0.119 |
| 265 | 3.922 | 1.927 | 0.078 |
| 266 | 3.962 | 0.927 | 0.037 |
| 266.927 | 4.000 | 0.000 | 0.000 |
| 268 | 4.043 | 1.073 | 0.043 |
| 269 | 4.083 | 2.073 | 0.082 |
| 270 | 4.122 | 3.073 | 0.122 |
| 271 | 4.161 | 4.073 | 0.161 |
| 272 | 4.200 | 5.073 | 0.200 |

This function indicates that a measurement error of 1 % in the flow value $\hat{q}$ translates into an error of around 3 % in the estimation value $\hat{l}$ of the contact point.

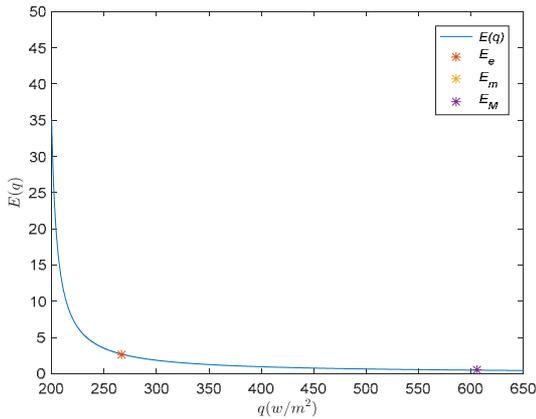

**Figure 5.** Elasticity for Example 2

### 5.3 Example 3

In this example we consider case where conductivity values are of the similar order. In particular, a Al-Mg bar is assumed to estimate the interface length length $l$ ($\kappa_A = 204 \text{Wm}^{-1}\,{}^0C^{-1} < \kappa_B = 156 \text{Wm}^{-1}\,{}^0C^{-1}$).

For $q = 474.475$ W.m$^{-2}$, the necessary and sufficient condition (31) are verified since that $q_m = 457.031$ W.m$^{-2}$ and $q_M = 503.289$ W. m$^{-2}$ and from equation (24) it follows that $l = 4m$.

Note that since the materials have similar conductivity values, the measurement of the heat flux should be more accurate in order to satisfy the condition (31).

The Tab. 4 contains the noisy value for flow measurements in a range around the actual value $q$, along the level noise $\varepsilon$ with the corresponding estimates $\hat{l}$ for $l$ and the estimated error bound $K$ defined in (35).

The estimate value in Tab. 4 are less accurate than those obtained in examples 1 and 2. This is due to the fact that the thermal conductivities of the materials are similar and by equation (35) this causes the error in the estimation to increase. For the same reason $K$ values are higher in this example compared to the ones for examples 1 and 2.

**Table 4.** Estimate of $l$ for Example 3

| $\hat{q}$ (W. m$^{-2}$) | $\hat{l}$ (m) | $\varepsilon$ (W. m$^{-2}$) | $K$(m) |
|---|---|---|---|
| 470 | 3.002 | 4.475 | 0.998 |
| 471 | 3.226 | 3.475 | 0.773 |
| 472 | 3.450 | 2.475 | 0.549 |
| 473 | 3.673 | 1.475 | 0.327 |
| 474 | 3.895 | 0.475 | 0.105 |
| 474.475 | 4.000 | 0.000 | 0.000 |
| 476 | 4.336 | 1.525 | 0.336 |
| 477 | 4.555 | 2.525 | 0.555 |
| 478 | 4.773 | 3.525 | 0.773 |
| 479 | 4.899 | 4.525 | 0.990 |

The estimate values in Tab. 4 are less accurate than those obtained in examples 1 and 2. This is due to the fact that the thermal conductivities of the materials are similar and by equation (35) this causes greater errors in the estimation. For the same reason the $K$ values (errors bounds) are higher in this example compared to the ones obtained for examples 1 and 2.

In Fig. 6 the elasticity is plotted for example 3. Notice that for this example $\kappa_A > \kappa_B$, and for this reason, elasticity is a positive function (see (40)), strictly decreasing (see (41)) having a vertical asymptote at $q = q_m$ (see (38)).

This function indicates that a measurement error of 1 % in the flow value $\hat{q}$ translates into an error higher than 25 % in the estimation value $\hat{l}$ of the contact point.

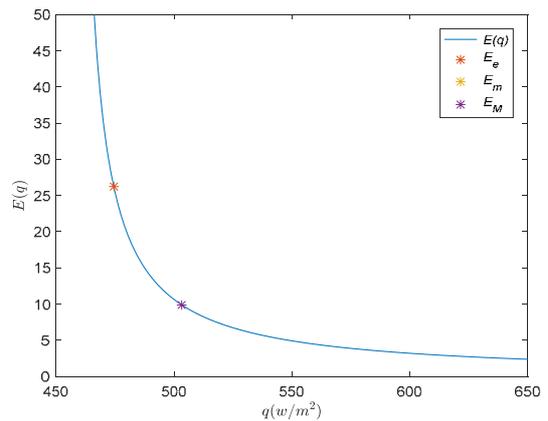

**Figure 6.** Elasticity for Example 3

The above examples show the different situations related to the conductivity values of the two materials. In all cases, the necessary and sufficient conditions (31) are verified. Although particular materials and a particular value for the interface location are considered, other cases have similar behavior.

### 6. CONCLUSIONS

A mathematical model for the heat transfer of a bar of negligible diameter and known length totally isolated on its lateral surface composed of two different, isotropic and homogeneous materials is studied. Appropriate boundary conditions are imposed and an analytical solution to the forward problem is found that turns out to be consistent with the case where no interface is present (a bar made of only one material).

Using a flow over-specified condition on the far right, a technique for estimating the contact point is proposed. Necessary and sufficient conditions are provided for the estimation of the parameter as well as a bound for the estimation error.

Using the elasticity function, the local influence of the dependency of the contact point on the flow is studied.

The numerical examples indicate that the approach introduced here is useful for determining the point of contact between the materials, but it is necessary that the flow be measured as accurately as possible in order to obtain a good estimation. The elasticity analysis indicates that the error estimates could become very large when the materials have similar thermal conductivities.

## ACKNOWLEDGMENT

This work has been partially supported by SOARD/AFOSR through grant FA9550-18-1-0523 for the first and second authors. This work has also been partially supported by European Union's Horizon 2020 Research and Innovation Program under the Marie Sklodowska-Curie Grant Agreement No. 823731 CONMECH and by the Project PIP No. 0275 from CONICET – Universidad Austral, Rosario, Argentina for the third author.

## REFERENCES

[1] Chung, D.D.L. (2001). Thermal interface materials. Journal of Materials Engineering and Performance, 10(1):56-59.
http://dx.doi.org/10.1361/105994901770345358

[2] Ma, H.Y., Zhu, P. X., Zhou, S. G., Xu, J., Huang, W. F., Yang, H.M., Chen, M.J. (2010). Preliminary research on Pb-Sn-Al laminated composite electrode materials Applied to zinc electrodeposition. Advanced Materials Research, 150: 303–308. http://dx.doi.org/10.4028/www.scientific.net/AMR.150-151.303.

[3] Cahill, D.G., Ford, W.K., Goodson, K.E., Mahan, G.D., Majumdar, A., Maris, H.J., Merlin, R., Phiipot, S.R. (2003). Nanoscale thermal transport. Journal of Applied Physics, 93(2): 793-818. http://dx.doi.org/10.1063/1.1524305.

[4] Clausing, A.M., Chao, B.T. (1965). Thermal contact resistance in a vacuum environment. Journal Heat Transfer, 87(2): 243-251. http://dx.doi.org/10.1115/1.3689082.

[5] Barturkin, V. (2005). Micro-satellites thermal control-concepts and components. Acta Astronautica, 56(1-2): 161-170.

[6] Shervedani, R.K., Isfahani, A.Z., Khodavisy, R., Hatefi-Mehrjardi, A. (2007). Electrochemical investigation of the anodic corrosion of Pb–Ca–Sn–Li grid alloy in H2SO4 solution. Journal of Power Sources, 164(2): 890-895.
http://dx.doi.org/10.1016/j.jpowsour.2006.10.105.

[7] Zhang, L., Yang, P., Chen, M., Liao, N. (2012). Numerical investigation on thermal properties at Cu–Al interface in micro/nano manufacturing. Applied Surface Science, 258(8): 3975–3979. http://dx.doi.org/10.1016/j.apsusc.2011.12.075.

[8] Yu, J., Yee, A.L., Schwall, R.E. (1992). Thermal conductance of Cu/Cu and Cu/Si interfaces from 85 K to 300 K. Cryogenics, 32(7): 610–615. http://dx.doi.org/10.1016/0011-2275(92)90291-H.

[9] Yang, P., Liao, N.B., Li, C., Shang, S.H. (2009). Multi-scale modeling and numerical analysis of thermal behavior of Cu-Al interface structure in micro/nano manufacturing. International Journal of Nonlinear Sciences and Numerical Simulation, 10(4): 483-491. http://dx.doi.org/10.1515/IJNSNS.2009.10.4.483.

[10] Prabhu, K.N., Kumar, S.T., Venkataraman, N. (2002). Heat transfer at the metal/substrate interface during solidification of Pb-Sn solder alloys. Journal of Materials Engineering and Performance, 11(3): 265–273.
http://dx.doi.org/10.1361/105994902770344051.

[11] Ramesh, G., Prabhu, K.N. (2012). Heat transfer at the casting/chill interface during solidification of commercially pure Zn and Zn base alloy (ZA8). International Journal of Cast Metals Research, 25(3): 160-164.
http://dx.doi.org/10.1179/1743133611Y.0000000026.

[12] Silva, J.N., Moutinho, D.J., Moreira, A.L., Ferreira, I.L., Rocha, O.L. (2011). Determination of heat transfer coefficients at metal–mold interface during horizontal unsteady-state directional solidification of Sn–Pb alloys. Materials Chemistry and Physics, 130(1-2): 179–185. http://dx.doi.org/10.1016/j.matchemphys.2011.06.032.

[13] Volz, S., Saulnier, J., Cheng, G., Beauchamp, P. (2000). Computation of thermal conductivity of Si/Ge superlattices by molecular dynamics techniques. Microelectronics Journal, 31 (9-10): 815-819. http://dx.doi.org/10.1016/S0026-2692(00)00064-1

[14] Ward, D.K., Curtin, W.A., Qi, Y. (2006). Aluminum–silicon interfaces and nanocomposites: A molecular dynamics study. Composites Science and Technology, 66(9): 1151-1161. http://dx.doi.org/10.1016/j.compscitech.2005.10.024

[15] Han, Z.H., Zhu, P.X., Guo, Y.Z., Zhou, S.G. (2015). Study on the interface and performance of Ti-Al laminated composite electrode materials. Journal of New Materials for Electrochemical Systems, 18(3): 159-164. https://doi.org/10.14447/jnmes.v18i3.363

[16] Yang, P., Liao, N. (2007). Surface sliding simulation in micro-gear train for adhesion problem and tribology design by using molecular dynamics model. Computational Materials Science, 38(4): 678–684. http://dx.doi.org/10.1016/j.commatsci.2006.06.004.

[17] Rogers, G.F.C. (1961). Heat transfer at the interface of dissimilar metals. International Journal of Heat and Mass Transfer, 2(1-2): 150–154. http://dx.doi.org/10.1016/0017-9310(61)90022-9.

[18] Cengel, Y., Heat and mass transfer. A practical approach. (2007). McGraw Hill, New York.

[19] Özisik, M.N., Heat conduction. (1993). John Wiley & sons, New York.

[20] Banks, H.T., Dediu, S., Ernstberger, S.L. (2007). Sensitivity functions and their uses in inverse problems. Journal Inverse and Ill-posed Problems, 15: 683-708. http://dx.doi.org/10.1515/jiip.2007.038.

[21] Sydsaeter, K., Hammond, P.J. (1995). Mathematics for economic analysis, Prentice Hall, New Jersey.

## NOMENCLATURE

| | |
|---|---|
| $a$ | variable assistant, °C |
| $b$ | variable assistant, °C.m$^{-1}$ |
| $c$ | variable assistant, °C |
| $d$ | variable assistant, °C.m$^{-1}$ |
| $E$ | elasticity |
| $F$ | heat source, °C |
| $h$ | convective coefficient, W. m$^{-2}$.°C$^{-1}$ |
| $K$ | bound for estimation error, m |
| $l$ | interface position, m |
| $\hat{l}$ | Interface position estimation, m |
| *L* | bar length, m |
| *q* | thermal flow, W.m$^{-2}$ |
| $\hat{q}$ | measured thermal flux, W.m$^{-2}$ |
| *S* | sensitivity, W$^{-1}$.m$^{3}$ |
| *Ta* | room temperature, °C |
| *u* | bar temperature, °C |
| *x* | special variable, m |

**Greek symbols**

| | |
|---|---|
| $\varepsilon$ | bound for flow measurement error, W.m$^{-2}$ |
| $\kappa$ | thermal conductivity, W.m$^{-1}$.°C$^{-1}$ |
| $\varsigma$ | variable assistant, W$^{2}$.m$^{-2}$.°C$^{-2}$ |

**Subscripts**

| | |
|---|---|
| A | regarding material A |
| B | regarding material B |
| e | regarding to the exact value |
| M | regarding the maximum |
| m | regarding the minimum |